\documentclass{amsart}
\newcommand{\modd}[1]{\, (\text{mod}\, {#1})}

\begin{document}
\parskip10pt
\parindent10pt
\baselineskip15pt

\newtheorem{theorem}{Theorem}
\newtheorem{lemma}{Lemma}

\title{Multiplicative Functions and A Taxonomy of Dirichlet Characters}

\author{P.D.T.A. Elliott}
\address{Department of Mathematics, University of Colorado Boulder, Boulder, Colorado 80309-0395 USA}
\email{pdtae@euclid.colorado.edu}

\author{Jonathan Kish}
\address{Department of Mathematics, University of Colorado Boulder, Boulder, Colorado 80309-0395 USA}
\email{jonathan.kish@colorado.edu}


\begin{abstract}
A method of estimating sums of multiplicative functions braided with Dirichlet characters is demonstrated, leading to a taxonomy of the characters for which such sums are large.
\end{abstract}

\maketitle

In this paper $g$ will denote a multiplicative arithmetic function with values in the complex unit disc.  In the seventh of a series of papers devoted to the study of such functions on arithmetic progression, c.f. Elliott \cite{elliottMFoAP1}, \cite{elliottMFoAP2}, \cite{elliottMFoAP3}, \cite{elliottMFoAP4}, \cite{elliottMFoAP5}, \cite{elliottMFoAP6}, \cite{elliottMFoAP7}, the first author showed that for  $\alpha > \frac{1}{4}$ there are non-principal Dirichlet characters $\chi_j \modd{D}$, $j = 1, \dots, k$ providing a representation
\begin{equation}
\sum\limits_{\substack{n \le x \\ n \equiv a \modd{D}}} g(n) = \frac{1}{\phi(D)} \sum\limits_{\substack{n \le x \\ (n,D)=1}} g(n) + \sum\limits_{j=1}^k \frac{\overline{\chi_j(a)}}{\phi(D)} \sum\limits_{n \le x} g(n)\chi_j(n) + O\left(\frac{x}{\phi(D)}\left(\frac{\log D}{\log x}\right)^{\alpha}\right) \notag
\end{equation}
uniformly for $(a,D)=1$, $2 \le D \le x$.  There are other uniformities too, for which we refer to the original papers.  Here $k$, the number of characters, depends at most upon $\alpha$.  It is further shown that with at most one exception the sums involving the characters $\chi_j$ fall within the bound on the error term.  One such exception is to be expected, since $g$ may itself be a character $\modd{D}$.\label{NEW.QUESTION.00}

Appealing, for example, to the Theorem in Elliott and Kish \cite{elliottandkishone2012}, the argument of Elliott \cite{elliottMFoAP6}, \cite{elliottMFoAP7}, yields a representation of the above type for each $\alpha > \frac{1}{2}$, and a variant of the argument allows any $\alpha < 1$ to be taken.  We do not further pursue this precision which is expected to be considered in a forthcoming book by A. Granville and K. Soundararajan.  The number of exceptional characters $\chi_j$ may be expected to increase as $\alpha$ increases but it will still depend at most upon $\alpha$.  The characters $\chi_j$ are otherwise unspecified.  The corresponding braided sums over $g\chi_j$ may no longer fall within the umbrella of the error term.

In the present paper we introduce a second method of estimating sums of multiplicative functions braided with Dirichlet characters.  Moreover, the method yields a taxonomy of the exceptional characters, classifying according to the size of a support set of the function, $g$, on the primes, i.e., a set of primes, $p$, outside of which $g(p)$ vanishes.  For applications to harmonic analysis it is important that the classification not depend upon the structure of the support set since the characterization of an extremal is often the aim of an application.

Functions of varying support that are restrictions of a given function $g$ are considered in Elliott \cite{elliottMFoAP7}; we shall argue from first principles.

The two methods may, with advantage, be combined.

In addition, we refine the dependence of the error term upon the modulus $D$.

A foundation is provided by


\begin{theorem}\label{thmelliottandkishtwo2012.1}
To each positive real $B$ there is a further real $c$ with the following property:

If $\chi$ is a Dirichlet character $\modd{D}$, $D \ge 2$, $8D^{-3B} \le \delta \le 1$, $t$, real, satisfies $|t| \le D^B$ and $S$ is a set of primes $p$ in the interval $(D,x]$ for which
\begin{equation}
\sum\limits_{p \in S} p^{-1}\left|1-\chi(p)p^{it}\right|^2 \le \delta L, \notag
\end{equation}
where $L = \sum_{D < p \le x} p^{-1}$, then
\begin{equation}
\text{either} \ \sum\limits_{p \in S}p^{-1} \le 4\delta^{1/3}L + c, \ \text{or the order of $\chi$ is less than} \ 2\delta^{-1/3}. \notag
\end{equation}
\end{theorem}

In particular, if $S$ contains a positive proportion of the primes in $(D,x]$ however small (in an obvious sense), then the order of the character $\chi$ is bounded.  If the character has order $r$ then the factorisation $1-z^r = (1-z)\left(1+z+\cdots + z^{r-1}\right)$ allows that
\begin{equation}
\sum\limits_{p \in S} \left|1-p^{irt}\right|^2p^{-1} \le r^2 \sum\limits_{p \in S} \left|1-\chi(p)p^{it}\right|^2 p^{-1} \le r^2 \delta L, \notag
\end{equation}
so that there is a control upon the size of $|t|$, too, c.f. Elliott \cite{elliott1988additivearithfuncs}, Lemmas 6, 7.

Moreover, characters of a fixed order may be interpreted as power-residue symbols and studied using reciprocity laws.

As in Elliott \cite{elliott2002millenium}, \cite{elliottMFoAP6}, \cite{elliottMFoAP7}, Elliott and Kish, loc. cit., the following result is vital.


\begin{lemma}\label{lemelliottandkishtwo2012.1}
Given $B>0$,
\begin{equation}
\text{Re} \sum\limits_{w < p \le y} \chi(p)p^{-s} \notag
\end{equation}
is bounded above in terms of $B$ alone, uniformly for $s = \sigma + it$, $\sigma = \text{Re}(s) \ge 1$, $|t| \le D^B$, $y \ge w \ge D$ and all non-principal characters $\chi \modd{D}$, $D \ge 1$.

If $|t|\log D \ge 1$ then it is similarly bounded for the principal character $\modd{D}$. 
\end{lemma}

\emph{Proof of Lemma \ref{lemelliottandkishtwo2012.1}}.  Two proofs via complex variables are given in Elliott \cite{elliottMFoAP6}, an elementary proof via Selberg's sieve is given in Elliott \cite{elliott2002millenium}.


\emph{Proof of Theorem \ref{thmelliottandkishtwo2012.1}}.  We begin with the Fej\'er kernel\label{NEW.QUESTION.03}
\begin{equation}
\frac{1}{N}\left(\frac{\sin \pi N\theta}{\sin \pi\theta}\right)^2 = \sum\limits_{m = -(N-1)}^{N-1} \left(1-\frac{|m|}{N}\right)e^{2\pi i m \theta}, \quad \theta \in \mathbb{R}. \notag
\end{equation}
It is useful to note that if $\| \theta \|$ denotes the distance of $\theta$ to a nearest integer, then
\begin{equation}
2\|\theta\| \le \left|\sin \pi \theta\right| \le \pi \|\theta\|. \notag
\end{equation}
If, moreover, $2N\|\theta\| \le 1$, then $\|N\theta\| = N\|\theta\|$ and the Fej\'er kernel is at least $4N\pi^{-2}$.  With $2\pi\gamma_p = \arg\chi(p) + t\log p$ we see that
\begin{align}
4N\pi^{-2} \sum\limits_{\substack{p \in S \\ 2N\|\gamma_p\| \le 1}} \frac{1}{p} & \le \sum\limits_{D < p \le x}\frac{1}{N}\left(\frac{\sin \pi N\gamma_p}{\sin\pi\gamma_p}\right)^2 \frac{1}{p} \notag \\
& = \sum\limits_{m=-N+1}^{N-1} \left(1-\frac{|m|}{N}\right)\text{Re} \sum\limits_{D < p \le x} \left(\chi(p)p^{it}\right)^m \frac{1}{p}. \notag
\end{align}
If the order of $\chi$ is at least $N$ and $N|t|\le D^{2B}$, then each of the innermost character sums with $m \ne 0$ is by Lemma \ref{lemelliottandkishtwo2012.1} bounded above in terms of $B$ alone.  Hence
\begin{equation}
\sum\limits_{\substack{p \in S \\ 2N\|\gamma_p\| \le 1}} \frac{1}{p} \le \frac{\pi^2 L}{4N} + c, \notag
\end{equation}
for a suitable constant c.

However, $\left|1-\chi(p)p^{it}\right| = \left|1-e^{2\pi i \gamma_p}\right| = 2\left|\sin\pi\gamma_p\right| \ge 4\|\gamma_p\|$ so that
\begin{equation}
\sum\limits_{\substack{p \in S \\ 2N\|\gamma_p\| > 1}} \frac{1}{p} < \frac{N^2}{4} \sum\limits_{p \in S} \left|1-\chi(p)p^{it}\right|^2 \frac{1}{p} \le \frac{N^2 \delta L}{4}. \notag
\end{equation}
Altogether,
\begin{equation}
\sum\limits_{p \in S} \frac{1}{p} \le \left(\frac{N^2\delta}{4} + \frac{\pi^2}{4N}\right)L + c. \notag
\end{equation}
Choosing, for simplicity of exposition, $N=[2\delta^{-1/3}]$, so that $N>2\delta^{-1/3}-1 \ge \delta^{-1/3}$, we obtain the upper bound $4\delta^{1/3}L +c$.

Note that if $\delta \ge 8D^{-3B}$, then $N|t| \le 2\delta^{-1/3}D^B \le D^{2B}$, as required.

Employing the lower bound $\sin \theta \ge \theta - \frac{\theta^3}{6}$, valid for $0 \le \theta \le \frac{\pi}{2}$, and the consequent lower bound
\begin{equation}
\frac{1}{r}\left(\frac{\sin \pi r\theta}{\sin \pi \theta}\right)^2 \ge r\left(1-\frac{(\pi r \|\theta\|)^2}{6}\right) \notag
\end{equation}
for the Fej\'er kernel, valid if $2r\|\theta\| \le 1$, the same argument delivers the following variant of Theorem~\ref{thmelliottandkishtwo2012.1}.


\begin{theorem}\label{thmelliottandkishtwo2012.2}
Under the hypotheses of Theorem \ref{thmelliottandkishtwo2012.1}, if $r^3\delta \le 1$ and $\chi$ has order at least $r$, then
\begin{equation}
\sum\limits_{p \in S} \frac{1}{p} \le \left(1 + (r^3\delta)^{1/2}\right)\frac{L}{r} + c, \notag
\end{equation}
for a suitable constant $c$.
\end{theorem}

If $D$ is a prime, $D \equiv 1 \modd{r}$ and $S$ is the set of primes in $(D,x]$ that are $r^{th}$-powers$\modd{D}$, then for every character of order $r$: $\sum_{p \in S} p^{-1} \left|1-\chi(p)\right|^2 = 0$ whilst, for a fixed $D$, Dirichlet's Theorem on primes in arithmetic progressions shows that $\sum_{p \in S}p^{-1} = r^{-1}L + O(1)$ as $x \to \infty$.  The upper bound in Theorem \ref{thmelliottandkishtwo2012.2} cannot be appreciably improved.

We return to the taxonomy.  In the next section $\varepsilon$, $0 < \varepsilon \le \frac{1}{4}$, is considered fixed, $2 \le D \le x^{3/4}$; $S$, the support of the function $g$ on the primes, is contained in $(D,x]$; $L$ is $\sum_{D < p \le x}p^{-1}$ as before and $\beta = L^{-1}\sum_{p \in S}p^{-1}$.  For background we recall the following refinement of a Theorem of Hal\'asz, a proof of which may be found in Elliott \cite{elliottMFoAP7}, Lemma 9.


\begin{lemma}\label{lemelliottandkishtwo2012.2}
Let $T \ge 2$.  Let $g$ be a multiplicative arithmetic function for which $g(p)$ vanishes on the primes outside the interval $(Y,x]$, $\frac{3}{2} \le Y \le x$.  Then
\begin{equation}
x^{-1}\log Y \sum\limits_{n \le x}g(n) \ll (m+1)e^{-m} + T^{-1/2}, \notag
\end{equation}
where
\begin{equation}\label{QUESTION.11}
m = \min\limits_{|t| \le T} \sum\limits_{Y < p \le x} p^{-1}\left(1-\text{Re}\, g(p)p^{it}\right), \notag
\end{equation}
the implied constant absolute.
\end{lemma}

\emph{The taxonomy}.

\emph{Case 1}.  Suppose that $\beta \ge \frac{1}{2} + \varepsilon$, and that for distinct characters $\chi_j \modd{D}$,
\begin{equation}
\sum\limits_{p \in S} p^{-1}\left|1-g(p)\chi_j(p)p^{-it_j}\right|^2 \le \frac{1}{4}\delta L, \quad |t_j| \le D^B, j = 1, 2. \notag
\end{equation}
An application of the Cauchy-Schwarz inequality shows that
\begin{equation}
\sum\limits_{p \in S} p^{-1}\left|1-\chi_1\overline{\chi}_2(p)p^{i(t_1-t_2)}\right|^2 \le \delta L, \notag
\end{equation}
and by Theorem \ref{thmelliottandkishtwo2012.2} with $r=2$, $\varepsilon \le (2\delta)^{\frac{1}{2}} + cL^{-1}$ for a constant $c$ that depends upon $B$ alone.  Note that if $8\delta > 1$ then $(2\delta)^{1/2} > \frac{1}{2}$ and this last inequality is already satisfied.

In particular the inequality
\begin{equation}
\sum\limits_{p \in S} \left(1-\text{Re}\, g(p)\chi(p)p^{it}\right) p^{-1} > \frac{1}{2}\varepsilon^2 L - c \notag
\end{equation}
holds uniformly for $|t| \le D^B$ and all save possibly one character $\chi \modd{D}$.

Adapting the methods in Elliott \cite{elliottMFoAP6}, \cite{elliottMFoAP7}, which includes an appeal to Lemma \ref{lemelliottandkishtwo2012.2}, we then arrive at the estimate
\begin{align}
\sum\limits_{\substack{n \le x \\ n \equiv a \modd{D}}} g(n) & = \frac{1}{\phi(D)} \sum\limits_{\substack{n \le x \\ (n,D)=1}} g(n) + \frac{\overline{\chi_1(a)}}{\phi(D)} \sum\limits_{n \le x} g(n)\chi_1(n) +  \notag \\
&  \hspace{1in} + O\left(\frac{x}{\phi(D)} \left(\frac{\log D}{\log x}\right)^{1-\beta+\frac{\varepsilon^2}{8}} \frac{\log(\log x/\log D)}{\log D}\right) \notag
\end{align}
uniformly for $(a,D)=1$, with a single possible exceptional character $\chi_1$ and with various attached uniformities that we do not pursue here.

\emph{Case 2}.  We weaken the bound on $\beta$ to $\beta \ge \frac{1}{3} + \varepsilon$.  Little changes save that with $r=3$ we may only conclude that  $\chi_1\overline{\chi}_2$ has order at most 2.  The upshot is an estimate
\begin{align}
\sum\limits_{\substack{n \le x \\ n \equiv a \modd{D}}} g(n) & = \frac{1}{\phi(D)} \sum\limits_{\substack{n \le x \\ (n,D)=1}} g(n) + \sum\limits_{\chi^2 = \chi_0} \frac{\overline{\chi}_1\chi(a)}{\phi(D)}\sum\limits_{n \le x} g(n)\chi_1\chi(n)+ \notag \\
& \hspace{1in} + O\left(\frac{x}{\phi(D)}\left(\frac{\log D}{\log x}\right)^{1 - \beta + \frac{\varepsilon^2}{12}}\frac{\log(\log x/\log D)}{\log D}\right) \notag
\end{align}
with uniformities as before.  Here $\chi_0$ denotes the principal character$\modd{D}$.

Note that if $g$ is real-valued then the hypothesis $\sum_{p \in S} p^{-1}\left|1-g(p)\chi(p)p^{it}\right|^2 \le \left(\frac{\delta}{4}\right)L$ guarantees at once that
\begin{equation}
\sum\limits_{p \in S}p^{-1} \left|1-\chi(p)^2 p^{2it}\right|^2 \le \delta L. \notag
\end{equation}
We may assume the exceptional character $\chi_1$ in the first case to be real, and may delete it altogether in the second case.

It is clear how this process continues.

In the general case the exponent of the error term may be given the value $1-\beta + \frac{1}{4}\min\left(\frac{1}{r^3},\frac{\varepsilon^2}{r}\right)$.  Denoting this exponent by $1-\beta + \omega_r$, the upshot for a multiplicative function $g$ with \emph{unconstrained} support $S$ whose restriction to the interval $(D,x]$ satisfies $\beta \ge 1/r + \varepsilon$, is a representation\label{QUESTION.14}
\begin{align}
\sum\limits_{\substack{n \le x \\ n \equiv a \modd{D}}} g(n) & = \frac{1}{\phi(D)}\sum\limits_{\substack{n \le x \\ (n,D)=1}} g(n) + \frac{1}{\phi(D)}\sum\limits_{\substack{\chi \ \text{has order} \\  2, \dots, r-1}} \frac{\overline{\chi_1\chi}(a)}{\phi(D)}\sum\limits_{n \le x} g(n)\chi_1\chi(n) \notag \\
& \hspace{1in} + O\left(\frac{1}{\phi(D)}\frac{x}{\log x}\prod\limits_{\substack{p \le x \\ p \in S}} \left(1+\frac{1}{p}\right)\left(\frac{\log D}{\log x}\right)^{\omega_r}\log\left(\frac{\log x}{\log D}\right)\right), \notag
\end{align}
with attached uniformities.

\emph{Remarks}.

1.  When $g$ is the M\"obius function, $\mu$, appropriate to the study of primes in arithmetic progression, in particular to a derivation of Linnik's theorem on the size of the least prime in such a progression, c.f. Elliott \cite{elliottMFoAP6}, \cite{elliott2002millenium}, $\beta = 1$ and there can be at most one exceptional character, that real.

2.  When considering multiplicative functions with varying support on the primes, the following result is helpful, c.f. Elliott \cite{elliottMFoAP6}, Lemma 7.


\begin{lemma}\label{lemelliottandkishtwo2012.3}
Let $0 < \varepsilon < 1$.  The inequality
\begin{equation}
\sum\limits_{j=1}^J \max\limits_{v-u \le H} \left|\sum\limits_{\substack{u < n \le v \\ (n,Q)=1}} a_n\chi_j(n)\right|^2 \ll \left(H\prod\limits_{\substack{p \mid Q \\ p \le H}} \left(1-\frac{1}{p}\right) + JH^\varepsilon D^{\frac{1}{2}}\log D\right)\sum\limits_{n=1}^\infty |a_n|^2 \notag
\end{equation}
where the $\chi_j$ are distinct Dirichlet characters $\modd{D}$, $D \ge 2$, $Q$ a positive integer, $H \ge 0$, holds for all square-summable complex numbers $a_n$, the implied constant depending at most upon $\varepsilon$.
\end{lemma}

There are several ways to establish this result, which is of Maximal Gap Large Sieve type.  An application of Cauchy's inequality shows that we may assume $D^\frac{1}{2}$ not to exceed $H$.  We may also include in $Q$ the prime-divisors of $D$.  Note that $\sum\limits_{\substack{p \mid Q \\ x < p \le x^2}} p^{-1} \le \sum\limits_{x < p \le x^2} p^{-1} \ll 1$ uniformly in $x \ge 1$.


\emph{Proof of Lemma \ref{lemelliottandkishtwo2012.3}}.  With $0 \le v_j-u_j \le H$, define
\begin{equation}
t_j(n) = \begin{cases} \chi_j(n) & \text{if} \ u_j < n \le v_j, \\ 0 & \text{otherwise}, \end{cases} \qquad j = 1, \dots, J. \notag
\end{equation}
For any real $\lambda_d$, $d\mid Q$, constrained by $\lambda_1 = 1$, the dual form
\begin{equation}
S = \sum\limits_{\substack{n \le x \\ (n,Q)=1}} \left|\sum\limits_{j=1}^J c_jt_j(n)\right|^2 \notag
\end{equation}
does not exceed
\begin{equation}
\sum\limits_{n \le x} \left(\sum\limits_{d \mid (n,Q)} \lambda_d \right)^2 \left|\sum\limits_{j=1}^J c_jt_j(n)\right|^2 \notag
\end{equation}
\begin{equation}
= \sum\limits_{d_i \mid Q} \lambda_{d_1}\lambda_{d_2} \sum\limits_{j,k=1}^J c_j \overline{c}_k \sum\limits_{ n \equiv 0 \modd{[d_1,d_2]}} t_j(n) \overline{t_k(n)}, \notag
\end{equation}
where $[d_1,d_2]$ denotes the least common multiple of $d_1$ and $d_2$.  For those terms with $j \ne k$, the innermost sum has the form
\begin{equation}
\chi_j\overline{\chi}_k\left([d_1,d_2]\right)\sum\limits_{m} \chi_j\overline{\chi}_k(m) \notag
\end{equation}
with the integers $m$ over an interval and is, by a classical result of P\'olya and Vinogradov, $O\left(D^\frac{1}{2}\log D\right).$  The corresponding contribution to $S$ is
\begin{equation}
\ll JD^\frac{1}{2}\log D \left(\sum\limits_{d \mid Q} |\lambda_d|\right)^2 \sum\limits_{j=1}^J |c_j|^2. \notag
\end{equation}
For those terms with $j=k$ we reform the square in the $\lambda_d$ to gain a contribution
\begin{equation}
\sum\limits_{j=1}^J |c_j|^2 \sum\limits_{n \le x} \left(\sum\limits_{d\mid (n,Q)} \lambda_d\right)^2 \left|t_j(n)\right|^2. \notag
\end{equation}
Since $\left|t_j(n)\right| \le 1$, the innersum over $n$ does not exceed
\begin{equation}
\sum\limits_{d_i \mid Q} \lambda_{d_1}\lambda_{d_2} \sum\limits_{\substack{u_j < n \le v_j+H \\ n \equiv 0 \modd{[d_1,d_2]}}} 1 = H\sum\limits_{d_i \mid Q} \lambda_{d_1}\lambda_{d_2}[d_1,d_2]^{-1} + O\left(\left(\sum\limits_{d\mid Q} |\lambda_d|\right)^2\right). \notag
\end{equation}
We may follow the standard appeal to the method of Selberg, c.f. Elliott \cite{Elliott1979}, Chapter 2, with $\lambda_d=0$ if $d>H^{\varepsilon/2}$ which, in particular, gives $|\lambda_d| \le 1$ for all remaining $\lambda_d$.

As a consequence,
\begin{equation}
S \ll \left(H\prod\limits_{\substack{p \mid Q \\ p \le H}} \left(1-\frac{1}{p}\right) + JH^\varepsilon D^\frac{1}{2}\log D\right)\sum\limits_{j=1}^J |c_j|^2. \notag
\end{equation}
Dualising gives the inequality of Lemma \ref{lemelliottandkishtwo2012.3}.

\emph{Error dependence upon the modulus}.

For certain applications it is helpful that the factor $\phi(D)^{-1}$ in the error terms may be replaced by $D^{-1}$.  Since the argument initially operates upon a multiplicative function $g$ with $g(p)$ supported on primes in an interval $(D^c,x]$, it is only at the end of the proof, when this restriction is removed, that more care needs to be taken.  It suffices to note that
\begin{equation}
\sum\limits_{\substack{n \le x, (n,D)=1 \\ p \mid n \implies p \le D^c}} \frac{|g(n)|}{n} \le \prod\limits_{\substack{p \le D^c \\ (p,D)=1}} \left(1-\frac{1}{p}\right)^{-1} \ll \frac{\phi(D)\log D}{D}, \notag
\end{equation}
in place of the bound $\ll \log D$ employed in Elliott \cite{elliottMFoAP6}, \cite{elliottMFoAP7}; and that the number of integers, not exceeding $x$, that lie in a reduced class$\modd{D}$ and have prime factors not exceeding $D^c$ is $\ll D^{-1}x \exp\left(-\log x/(25\log D)\right)$, uniformly for $2 \le D \le x$.  This last may be effected by appeal to the general probability model in Elliott \cite{Elliott1979}, Chapter 3, in particular Lemma 3.5, in conjunction with Elliott \cite{Elliott1979}, Lemma 3.3, on the same lines as the argument in Elliott \cite{elliottMFoAP6}, Lemma 13.

An immediately accessible example, c.f. Elliott \cite{elliottMFoAP7}, is that given $\alpha > \frac{1}{4}$, $g$ multiplicative with values in the complex unit disc, for each modulus $D$, $2 \le D \le x$, there is a Dirichlet character $\chi_1 \modd{D}$ such that
\begin{equation}
\sum\limits_{\substack{n \le x \\ n \equiv h \modd{D}}} g(n) = \frac{1}{\phi(D)}\sum\limits_{\substack{n \le x \\ (n,D)=1}} g(n) + \frac{\overline{\chi_1}(h)}{\phi(D)}\sum\limits_{n \le x} g(n)\chi_1(n) + O\left(\frac{x}{D}\left(\frac{\log D}{\log x}\right)^\alpha\right), \notag
\end{equation}
uniformly for $(h,D)=1$, the implied constant depending at most upon $\alpha$.

The following result, sufficient but more modest, may be obtained without appeal to the construction of a probability space.


\begin{lemma}\label{lemelliottandkishtwo2012.4}
Let $c>0$.  For any positive integer $k$, the number of integers $n$ in the interval $x^{1/2}< n\le x$ made up of primes not exceeding $D^c$ and which satisfy $n \equiv a \modd{D}$ is $\ll D^{-1}x(\log D/\log x)^k$ uniformly in $D \le x$ and $a$.
\end{lemma}


\emph{Proof of Lemma \ref{lemelliottandkishtwo2012.4}}.  Let $(a,D)=1$, $D \le x^\frac{1}{2}$.  Suppose that the integers $n$ are squarefree.  Then their number does not exceed
\begin{equation}
\sum\limits_{\substack{n \le x \\ n \equiv a \modd{D} \\ p\mid n \implies p \le D^c}} \mu^2(n)\left(\frac{\log n}{\log x^\frac{1}{2}}\right)^k \notag
\end{equation}
\begin{equation}
=\left(\frac{1}{2}\log x\right)^{-k}\sum\limits_{p_j \le D^c}\log p_1 \cdots \log p_k \sum\limits_{\substack{m \le x[p_1,\dots,p_k]^{-1} \\ m[p_1,\dots,p_k]\equiv a \modd{D}}} 1. \notag
\end{equation}
Here $[p_1,\dots,p_k] \le D^{kc}$ and if $D^{kc} \le x^\frac{1}{4}$, say, then the innersum is $\ll D^{-1}x[p_1,\dots,p_k]^{-1}$.

The terms with the $p_j$ distinct contribute
\begin{equation}
\ll \frac{x}{D(\log x)^k}\left(\sum\limits_{p \le D^c} \frac{\log p}{p}\right)^k \ll \frac{x}{D}\left(\frac{\log D}{\log x}\right)^k. \notag
\end{equation}

Those terms with two $p_j$ equal contribute
\begin{equation}
\ll \frac{x}{D(\log x)^k} \left(\sum\limits_{p\le D^c} \frac{\log p}{p}\right)^{k-2} \cdot \sum\limits_{p \le D^c} \frac{(\log p)^2}{p}, \notag
\end{equation}
which has a similar bound; and so on.

If $D>x^\frac{1}{2}$ or $D^{kc} \ge x^\frac{1}{4}$ then the desired bound is evident, since we may then ignore the condition that the prime factors of $n$ do not exceed $D$.

More generally, any positive integer has a representation $rm^2$ with $r$ squarefree and we may estimate the number of such integers not exceeding $x$ which have only prime factors not exceeding $D^c$, and which satisfy $rm^2 \equiv a \modd{D}$.  For $D \le x^\frac{1}{8}$ those integers with $m \le x^\frac{1}{4}$ are, by our above argument,
\begin{equation}
\ll \sum\limits_{m \le x^\frac{1}{4}} \frac{x}{Dm^2}\left(\frac{\log D}{\log x}\right)^k \ll \frac{x}{D}\left(\frac{\log D}{\log x}\right)^k \notag
\end{equation}
in number.  Omitting side conditions, those with $m>x^\frac{1}{4}$ are
\begin{equation}
\ll \sum\limits_{m > x^\frac{1}{4}} m^{-2}x \ll x^\frac{3}{4} \ll D^{-1}x^\frac{7}{8} \notag
\end{equation}
in number.  Once again, the condition $D \le x^\frac{1}{8}$ may be removed in favour of $D \le x$.

If now $(a,D)=t$, then the integers $n$ in the Lemma have the form $tw$ where $w \le x/t$, $w$ has only prime factors not exceeding $D^c$ and $w \equiv t^{-1}a \modd{t^{-1}D}$.  For $D \le x^\frac{1}{2}$ their number is
\begin{align}
& \ll \frac{t^{-1}x}{t^{-1}D}\left(\frac{\log t^{-1}D}{\log t^{-1}x}\right)^k \notag \\
& \ll \frac{x}{D}\left(\frac{\log D}{\log x}\right)^k, \notag
\end{align}
and once again the condition $D \le x^\frac{1}{2}$ may be relaxed to $D \le x$.

\emph{Concluding remark}.

In the foregoing account the function $g$ is compared to the function that is identically 1 on the primes; the corresponding Dirichlet series $\sum_{n=1}^\infty g(n)n^{-s}$ is compared to $\zeta(\sigma)$, where $\sigma = \text{Re}(s)>1$.  One may also compare to $\sum_{n=1}^\infty \left|g(n)\right|n^{-\sigma}$.  The basic assumption of Theorem \ref{thmelliottandkishtwo2012.1} is replaced by
\begin{equation}
\sum\limits_{D<p \le x}p^{-1}\left|g(p)\right|\cdot \left|1-\chi(p)p^{it}\right|^2 \le \delta L, \notag
\end{equation}
the initial hypothesis of the taxonomy by
\begin{equation}
\sum\limits_{\substack{D<p \le x \\ g(p) \ne 0}} \left|g(p)\right|\cdot \left|1-\left|g(p)\right|^{-1}g(p)p^{it_j}\right|^2 \le \frac{\delta}{4}L, \notag
\end{equation}
and so on.  The taxonomy is then in terms of the size of $\sum\limits_{D < p \le x} p^{-1}\left|g(p)\right|$.

\bibliographystyle{amsplain}
\bibliography{MathBib}

\end{document}